\documentclass{article}
\usepackage{amsmath}
  \usepackage{paralist}
  \usepackage{graphics} 
  \usepackage{epsfig} 
 \usepackage[colorlinks=true]{hyperref}
\hypersetup{urlcolor=blue, citecolor=red}

\newtheorem{theorem}{Theorem}[section]

\newtheorem{lemma}[theorem]{Lemma}

\newtheorem{definition}[theorem]{Definition}
\newtheorem{remark}{Remark}
\newtheorem{assumption}{Assumption}
\newtheorem{notation}{Notation}

\usepackage{amssymb}
\usepackage{amsfonts}
\usepackage{verbatim}
\usepackage{float}
\usepackage{subfigure}
\usepackage{graphicx}
\def\Om{\Omega}
\def\RR{{\mathbb{R}}}
\def\NN{{\mathbb{N}}}
\def\Lr{\mathcal{L}}

\def\OmLr{\Omega_\mathcal{L}}

\def\strain{\mathcal{E}^\star}
\def\nl{\newline}

\begin{document}

\title{The non-Riemannian dislocated crystal: a tribute to Ekkehart Kr\"oner (1919-2000)}

\author {\scshape Nicolas Van Goethem\footnote{Email: vangoeth@ptmat.fc.ul.pt. Work supported by Funda\c{c}\~{a}o para a Ci\^{e}ncia e a Tecnologia (Ci\^{e}ncia 2007 \& FCT Project: PTDC/EME-PME/108751/2008)}\\
Universidade de Lisboa\\Faculdade de Ci\^{e}ncias, Departamento de Matem\'{a}tica\\
Centro de Matem\'{a}tica e Aplica\c{c}\~{o}es Fundamentais,\\
Av. Prof. Gama Pinto 2, 1649-003 Lisboa, Portugal}
\date{}
\maketitle

\begin{abstract}
This expository paper is a tribute to Ekkehart Kr\"oner's results on the intrinsic non-Riemannian geometrical nature of a single 
crystal filled with point and/or line defects. A new perspective on this old theory is proposed, intended to contribute to the debate around 
the still open Kr\"oner's question: ``what are the dynamical variables of our theory?''
\end{abstract}

\section{Introduction}
In the field of solid state physics, in particular physics of defects, the legacy of Ekkehart Kr\"oner who died ten years ago at 
the age of $81$, is invaluable. He has been actively publishing for $50$ years, mostly as a single author, on the physical 
understanding of defective solids, but also on their mathematical structure. One could make a distinction between a first series of paper 
\cite{KR55}-\cite{KR80} where he constructs an original approach to understand dislocations, and a later series \cite{KR90}-\cite{KR01} where he raises questions, while reporting new knowledge in the field.

Most of the theory can be found in the course \cite{KR80} but since Kr\"oner also distilled many comments, ideas, and computations along other 
publications, the idea of writing the present tribute grew up. It is especially intended to commemorate the $10^{th}$ 
anniversary of his death, in order, 
not to recall (because the author has no privileged relationship with Kr\"oner to do so), but to enlighten Kr\"oner's ideas and show how they 
are found rich enough by the author to be diffused, revisited and emphasized today. 

It should be pointed out that Anthony \cite{ANT70a,ANT70b} is one of Kr\"oner's direct 
students who also greatly contributed to understand defect lines (in particular, disclinations). Since then, many contributions to the field 
(nonlinear dislocations, dislocation motion, thermodynamic of defective crystals etc) have appeared, but surprisingly enough, few only cite Kr\"oner. 
This is probably due to the lack of real school following him, but also due to scientific reasons: indeed, Kr\"oner's theory is formulated in 
physical terms, but appeals to complex mathematical concepts, the combination of which is only rarely seen in the literature. 
It should be emphasized that de Le\'{o}n, Epstein, Lazar, Maugin and co-authors \cite{LEONEP2,EPMAUG,BUCAEP,LAZMAUG} (cf. the well-documented survey 
\cite{Maugin2003} and the references therein) have produced 
significant results not only by following, but especially by completing the ideas of Kr\"oner.

So, the present paper is intended to (i) collect and show Kr\"oner's results in the light of a new presentation, (ii) describe the 
non-Riemannian crystal and show how it can help to select appropriate deformation and
internal (thermodynamic) variables, 
(iii) participate to the debate around Kr\"oner's question: ``what are the dynamical variables of our theory?'' \cite{KR95}

It will be especially stressed that the crystal geometry and the physical laws governing defects are inseparable, as is the case in 
the  Einstein's General Theory of Relativity. However, we entirely agree with Noll when he writes \cite{NOLL67} that ``the geometry [must be] the natural
outcome, not the first assumption, of the theory'' (i.e., as in the \textit{Continuous Distribution of Dislocation} (CDD)
theory of Bilby at al. \cite{BILBY}). Many geometrical tools and mathematical theory required for a rigourous description of
the dislocated crystal geometry can be found in the landmark papers by Noll \cite{NOLL67} and Wang \cite{WANG67}, while also
pointing out a recent book on 
Continuum Mechanics in that spirit \cite{EP2010}. The approach followed here and detailled in \cite{VGD2009}-\cite{VG2011} is nonetheless distinct 
from the CDD theory.

Single crystals growing from the melt are considered where high temperature gradients are unavoidable and hence where point defects are present \cite{VGetal}. 
Moreover, since there are no
internal  boundaries, the defect lines can take in principle any orientation while forming either loops or lines ending at the crystal 
boundary. However, for the purpose of simplicity in the exposition of the theory, we will consider a tridimensional crystal filled with a network of rectilinear parallel
disclinations and/or dislocations. 

Particular to the chosen approach is the distinction between scales, where the macroscale is recovered from the mesoscale by a 
homogenization process: the singularities (i.e., the defect lines) have been 
erased and hence the density of defects (dislocations and/or disclinations) are recovered by means of 
smooth fields which we will show responsible for 
curvature and torsion of the crystal intrisic geometry. Also, the density of point defects will show responsible for 
the appearance of non-metric terms. In this
approach, only objective fields are considered 
to describe defective matter: they are defined across scales although their physical meaning might differ. Moreover, no
elasto-plastic decomposition and no prescription of any reference configuration are required, and there is no assumption of 
static equilibrium (vanishing stress divergence).

\section{Preliminar results at the continuum scale}\label{meso}

\begin{notation}
In this paper, a scalar, vector  or tensor of any order are not typographically distinct symbols in the text. The tensor order is specified 
when equations are written, since in this case only, the vector $v$ is written as $v_i$ (with one index), and the tensor $U$ as $U_{ij\cdots}$ with a 
number of indices corresponding to its order.
\end{notation}

The present section focuses on the mesoscopic scale, where dislocations and disclinations are lines and whose 
characteristic length is some average distance between neighboring defects. The remaining of the medium is a 
continuum governed by linear elasticity. At time $t$, the body is referred to as ${\mathcal R}^{\star}(t)$ to represent any random 
sample corresponding to a given crystal growth experiment.  In the crystal domain $\Om$, the meso-scale physics will then be represented by a nowhere dense set of defect lines 
which in $2D$ are parallel to each other.

\begin{definition}[$2D$ mesoscopic defect lines]\label{lines2D}
At the meso-scale, a $2D$ set of dislocations and/or disclinations $\Lr\subset\Om$ is a closed set of $\Om$ 
(this meaning the intersection with $\Om$ of a closed set of $\RR^3$) formed by a countable union of parallel lines 
$L^{(i)}, i\in\mathcal{I}\subset\NN$, whose adherence is itself a countable union of lines  and where the linear 
elastic strain is singular. In the sequel, these lines will be assumed as parallel to the $z$-axis. 
\end{definition}

Since accumulation points (to be understood as clusters of parallel lines) might appear, the scale of matter
 description of this section is named continuum scale.

\subsection{Objective internal fields for the model description}

The present mesoscopic theory is developed from the sole linear elastic strain, which itself is
defined from the stress field (although the stress-strain relationship is not used in the sequel) and therefore is 
an objective internal field. 

\begin{assumption}[$2D$ mesoscopic elastic strain]\label{asstrain}
The linear strain $\mathcal{E}^\star$ is a given symmetric $L^1_{loc}(\Om)$ tensor such that $\partial_z\mathcal{E}^\star
=0$. Moreover, $\mathcal{E}^\star$ is assumed as compatible on $\OmLr:=\Om\setminus\Lr$ in the sense that the incompatibility tensor defined by
\begin{eqnarray}
&&\hspace{-150pt}\mbox{\scriptsize{INCOMPATIBILITY:}}\hspace{45pt}\eta^\star_{kl}:=\epsilon_{kpm}\epsilon_{lqn}\partial_p\partial_q\mathcal{E}^\star_{mn},\label{eta}
\end{eqnarray}
where derivation is intended in the distribution sense, vanishes everywhere on $\OmLr$. 
\end{assumption}

In the following definition generalizing the concept of rotation and displacement gradients to dislocated media, 
the strain is considered as a distribution on $\Omega$
 (i.e. as acting on $\mathcal{C}^1_c$ test-functions with compact support).
\begin{definition}[Frank and Burgers tensors]\label{FBtensors}
\begin{eqnarray}
&&\hspace{-118pt}\mbox{\scriptsize{FRANK TENSOR:}}
\hspace{45pt}\overline\partial_m\omega_k^\star:=\epsilon_{kpq}\partial_p\mathcal{E}_{qm}^\star\label{delta_m_a}\\
&&\hspace{-118pt}\mbox{\scriptsize{BURGERS TENSOR:}}
\hspace{41pt}\overline\partial_l b^\star_k:=\mathcal{E}^\star_{kl}+\epsilon_{kpq}(x_p-x_{0p})\overline\partial_l
\omega^\star_q,\label{delta_lb_i}
\end{eqnarray}
where $x_0$ is a point where displacement and rotation are given.
\end{definition}
Line integration of the Frank and Burgers tensors in $\OmLr$ (i.e., outside the defect set) provides the 
multivalued rotation and Burgers vector fields $\omega^\star$ and $b^\star$. These properties are summarized in the following 
theorem, whose proof is classical.

\begin{theorem}[Multiple-valued displacement field]\label{MVdispl}
From a symmetric smooth linear strain $\mathcal{E}^{\star}_{ij}$ on $\Om_\Lr$ and a point $x_0$ where displacement 
and rotation are given, a multivalued displacement field $u^\star_i$ can be constructed on $\Om_\Lr$ such that the symmetric part of the distortion $\partial_ju^\star_i$ is the single-valued strain tensor 
$\mathcal{E}^\star_{ij}$ while its skew-symmetric part is the multivalued rotation tensor
 $\omega^\star_{ij}:=-\epsilon_{ijk}\omega^\star_k$. Moreover, inside $\OmLr$ the gradient $\partial_j$ of 
the rotation and Burgers fields $\omega^\star_k$ and $b^\star_k=u^\star_k-\epsilon_{klm}\omega^\star_l(x_m-x_{0m})$ 
coincides with the Frank and Burgers tensors.
\end{theorem}

From this result, the Frank and Burgers vectors can be defined as invariants of any isolated defect line $L^{(i)}$ of 
$\Lr$.
\begin{definition}[Frank and Burgers vectors]\label{Burgers}
The Frank vector of the isolated defect line $L^{(i)}$ is the invariant
\begin{eqnarray}
\Om^{\star(i)}_k:=[\omega^\star_k]^{(i)},\label{frank}
\end{eqnarray}
while its Burgers vector is the invariant
\begin{eqnarray}
B^{\star(i)}_k:=[b^\star_k]^{(i)}=[u_k^\star]^{(i)}(x)-\epsilon_{klm}\Om^{\star(i)}_l(x_m-x_{0m}),\label{burgers}
\end{eqnarray}
with $[\omega^\star_k]^{(i)},[b^\star_k]^{(i)}$ and $[u_k^\star]^{(i)}$ denoting the jumps of $\omega^\star_k,b^\star_k$ 
and $u_k^\star$ around $L^{(i)}$.
\end{definition}
The three types of $2D$ defects are the screw and edge dislocation, and the wedge disclination.
As an example, the distributional strain and Frank tensor of an isolated screw dislocation (see \cite{VGD2009} for the 
other two) is given by the following 
$L^1_{loc}(\Om)$ symmetric tensor and Radon measure \cite{AFP2000}:
\begin{eqnarray}
[\mathcal{E}^\star_{ij}]&=&\frac{-B^\star_z}{4\pi r^2}\left[ \begin{array}{ccc} 0& 0 & y\\ 0 & 0 & -x \\ y & -x & 0
\end{array} \right]\nonumber\\
{[\overline\partial_m\omega^\star_k]}&=&
\frac{-B^\star_z}{4\pi r^{2}}\left[\begin{array}{ccc} \cos2\theta & \sin2\theta & 0 \\ \sin2\theta &-\cos2\theta   & 0 \\ 0 & 0 & 0
\end{array}\right]+
\frac{B^\star_z}{4}\left[\begin{array}{ccc} -\delta_L & 0 & 0 \\ 0 &-\delta_L   & 0 \\ 0 & 0 & 2\delta_L
\end{array}\right].\nonumber
\end{eqnarray}
Consider a set of countable lines $\Lr$ and remark that 
the present distributional approach is subtle in the sense that the physical condition that $\displaystyle\sum_{L^{(i)}\in\Lr}|B^{\star(i)}_z|$ 
be bounded is needed
in order for $\displaystyle\sum_{L^{(i)}\in\Lr}\overline\partial_m\omega^{\star(i)}_k$ to still be a Radon measure \cite{VGD2009}.

Besides their relationship with the multivalued rotation, Burgers and displacement fields, the Frank and Burgers 
tensors can be directly related to the strain incompatibility by use of (\ref{eta}), (\ref{delta_m_a}) \& 
(\ref{delta_lb_i}).
\begin{theorem}
The distributional curls of the Frank and Burgers tensors are
\begin{eqnarray}
\epsilon_{ilj}\partial_l\overline\partial_j\omega^\star_k&=&\eta^\star_{ik}\label{sup1}\\
\epsilon_{ilj}\partial_l\overline\partial_jb^\star_k&=&\epsilon_{kpq}(x_p-x_{0p})\eta^\star_{iq},\label{sup2}
\end{eqnarray} 
with $\eta^\star_{ik}$ the incompatibility tensor.
\end{theorem}

From this theorem it results that single-valued rotation and Burgers fields $\omega^\star$ and $b^\star$ can 
be integrated on $\Om$ if the incompatibility tensor vanishes.
\nl

To complete the model two other objective internal fields are introduced: the dislocation 
and disclination densities.
\begin{definition}[Defect densities]
\begin{eqnarray}
&&\hspace{-58pt}\mbox{\scriptsize{DISCLINATION DENSITY:}}\hspace{35pt}\Theta^\star_{ij}:=
\sum_{k\in\mathcal{I}\subset\NN}\Om^{\star (k)}_j\tau_i^{(k)}\delta_{L^{(k)}}\ (i,j=1\cdots 3)\label{disclindens1}\\
&&\hspace{-58pt}\mbox{\scriptsize{DISLOCATION DENSITY:}}\hspace{38pt}\Lambda^\star_{ij}:=
\displaystyle\sum_{k\in\mathcal{I}\subset\NN}B^{\star (k)}_j\tau_i^{(k)}\delta_{L^{(k)}}\ (i,j=1\cdots 3),\label{dislocdens2}
\label{dislocdens4}
\end{eqnarray}
where $\delta_{L^{(k)}}$ is used to represent the one-dimensional Hausdorff measure density \cite{AFP2000} 
concentrated on the rectifiable arc $L^{(k)}$ with the tangent vector $\tau_i^{(k)}$ defined almost everywhere 
on $L^{(k)}$, while $\Om^{\star (k)}_j$ and $B^{\star (k)}_j$ denote the Frank and Burgers vectors of $L^{(k)}$, 
respectively. 
\end{definition}

\subsection{Kr\"oner's formula}

In this paper, only a simplified $2D$ mesoscopic distribution of defects in a tridimensional crystal is considered. Accordingly, the vectors $\eta^\star_k,\Theta^\star_k$ and $\Lambda^\star_k$ denote the tensor components $\eta^\star_{zk},\Theta^\star_{zk}$ 
and $\Lambda^\star_{zk}$. 
Greek indices are used to denote the values $1,2$ (instead of the Latin indices used in $3D$ to denote the values 
$1,2$ or $3$). Moreover, $\epsilon_{\alpha\beta}$ denotes the permutation symbol $\epsilon_{z\alpha\beta}$.

The contortion as introduced by Kondo \cite{KONDO52}, Nye \cite{NYE53} and Bilby et al. \cite{BILBY} will show a crucial defect density tensor. 
Kr\"oner \cite{KR80} understood the importance of this object in terms of modelling.
\begin{definition}[$2D$ mesoscopic contortion]
\begin{eqnarray}
&&\hspace{-72pt}\mbox{\scriptsize{CONTORTION:}}\hspace{87pt}\kappa_{ij}^\star:=\delta_{iz}\alpha^\star_j
-\frac{1}{2}\alpha^\star_z\delta_{ij}\quad(i,j=1\cdots 3),\label{KR2} 
\end{eqnarray} 
where
\begin{eqnarray}
&&\hspace{68pt}\alpha^\star_j:=\Lambda^\star_j-\delta_{j\alpha}\epsilon_{\alpha\beta}\Theta^\star_z(x_\beta-x_{0\beta}).
\label{KR3}
\end{eqnarray} 
\end{definition}

Among several equivalent formulations, the formula relating strain incompatibility to the defect densities has been 
proposed in full generality by Kr\"oner \cite{KR80}, and proved for a countable set of $2D$ lines 
(this meaning, replacing subscript $i$ by $z$ in the formula) by Van Goethem \& Dupret \cite{VGD2009}:
\begin{eqnarray}
&&\hspace{-124pt}\mbox{\scriptsize{KR\"ONER'S FORMULA IN $2D$:}}\hspace{50pt}\eta_k^\star=\Theta_k^\star
+\epsilon_{\alpha\beta}\partial_\alpha\kappa^\star_{k\beta}.\label{etak}
\end{eqnarray}
For the expression of incompatibility for a set of skew isolated $3D$ lines, we refer to \cite{VG2010}, while general 
$3D$ results can be found in \cite{VG2011}.
\section{Preliminar results at the macroscopic scale}

Following Kondo \cite{KONDO55}, by calling a crystal ``perfect'', it is meant that the atoms form, in its stress-free 
configuration, a regular pattern proper to the prescribed nature of the matter. However, no real crystal is perfect, but 
rather filled with point and line defects which interact mutually. Each defect type is responsible for a 
particular geometric property, as will be described in this paper. In order to reach this crystal, we first need 
to provide a way from passing from the above scale to a scale where the fields have been smoothed.

\subsection{Homogenization}\label{sec_homo}
Homogenization is obtained from the continuum scale by a limit procedure which will not be detailled here (cf. \cite{VG2011}), 
but whose 
effect is to erase the singularities (isolated ones or those resulting from accumulation) and hence to 
provide a smooth macroscopic crystal. Basically we postulate the following limits:
\begin{equation}
 \Theta^\star,\Lambda^\star, \strain\rightarrow \Theta,\Lambda, 
\mathcal{E},\label{homo}
\end{equation} 
where $\Theta,\Lambda,\mathcal{E}$ belong to $\mathcal{C}^\infty(\Om)$ and where convergence is intended
in the sense of measures \cite{AFP2000}. The tensor $\mathcal{E}$ will be called macroscopic strain without 
claiming however that $\mathcal{E}$ is the elastic strain (i.e. linearly related to the macroscopic stress $\sigma$). 
As a consequence of law (\ref{homo}) we directly obtain from (\ref{eta}), (\ref{etak}) and straightforward 
distribution properties:
\begin{eqnarray}
&&\hspace{-42pt}\mbox{\scriptsize{MACROSCOPIC KR\"ONER'S FORMULA:}}\hspace{35pt}\eta_k=\epsilon_{\alpha\beta}\partial_\alpha
\overline\partial_\beta
\omega_k=\Theta_k
+\epsilon_{\alpha\beta}\partial_\alpha\kappa_{k\beta},\label{etakmacro}
\end{eqnarray}
where by (\ref{delta_m_a}) and (\ref{KR2}), (\ref{KR3}),
\begin{eqnarray}
\mbox{\scriptsize{MACROSCOPIC FRANK TENSOR:}}\hspace{12pt}\overline\partial_m\omega_k:=\epsilon_{kpq}\partial_p\mathcal{E}_{qm}\hspace{115pt}\label{frankmacrotens}\\
\hspace{-25pt}\mbox{\scriptsize{MACROSCOPIC CONTORTION:}}\hspace{25pt}\kappa_{ij}=\delta_{iz}\left(\Lambda^\star_j
-\delta_{j\alpha}\epsilon_{\alpha\beta}\Theta_z(x_\beta-x_{0\beta})\right)
-\frac{1}{2}\Lambda_z\delta_{ij}.
\label{dislocdens2d}
\end{eqnarray}

\subsection{External and internal observers}
The \textit{external observer} analyzes the crystal actual configuration $\mathcal{R}(t)$ with the Euclidian metric $g^{ext}_{ij}=\delta_{ij}$. 
The internal observer, in turn, can only count atom steps while moving in $\mathcal{R}(t)$, and parallelly transport a vector along crystallographic 
 lines. According to Kr\"{o}ner \cite{KR90}: ``in our universe we are internal observers who do not possess the ability to realize external 
actions on the universe, if there are such actions at all. Here we think of the possibility that the universe could be 
deformed from outside by higher beings. A crystal, on the other hand, is an object which certainly can deform from 
outside. We can also see the amount of deformation just by looking inside it, e.g., by means of an electron microscope. 
Imagine some crystal being who has just the ability to recognize crystallographic directions and to count lattice steps 
along them. 
Such an \textit{internal observer} will not realize deformations from outside, and therefore will be in a situation 
analogous to that of the physicist exploring the world. The physicist clearly has the status of an internal observer.''

\section{The macroscopic crystal}

At time $t$, the defective crystal is a tridimensional body denoted by $\mathcal{R}(t)$. The crystal 
defectiveness is not countable anymore, as was the case in Section \ref{meso}, since the fields have been smoothed 
by homogenization. However, defectiveness is recovered
by the natural embedding of the crystal into a specific geometry which will be described in the two following
sections.

\subsection{Macroscopic strain and contortion as key physical fields}

The macroscopic strain $\mathcal{E}$ and contortion $\kappa$ have been defined by homogenization in Section 
\ref{sec_homo}. It turns out that the relevant physical fields are not the Frank and Burgers tensors but their completed counterparts 
\cite{VGD2009}:

\begin{definition}\label{defmeso}
\begin{eqnarray}
\hspace{-45pt}\mbox{\scriptsize{COMPLETED FRANK TENSOR }}\hspace{38pt}\eth_j\omega_k
&:=&\overline\partial_j\omega_k-\kappa_{kj}\label{disclindens}\\
\hspace{-45pt}\mbox{\scriptsize{COMPLETED BURGERS TENSOR}}\hspace{33pt}\eth_j b_k
&:=&\mathcal{E}_{kj}
+\epsilon_{kpq}(x_p-x_{0p})\eth_j\omega_q.\label{dislocdens}
\end{eqnarray}
\end{definition}
The following result is a direct consequence of its mesoscopic counterpart:
\begin{theorem}\label{nouveautenseur1}
\begin{eqnarray}
\hspace{-42pt}\mbox{\scriptsize{MACROSCOPIC DISCLINATION DENSITY}}\hspace{65pt} 
\Theta_{ik}=\epsilon_{ilj}\partial_l\eth_j\omega_k&&\label{disclindens2}\\
\hspace{-34pt}\mbox{\scriptsize{MACROSCOPIC DISLOCATION DENSITY}}\hspace{68pt}
 \Lambda_{ik}=\epsilon_{ilj}\partial_l\eth_j b_k.&&\label{dislocdens2c}
\end{eqnarray}
\end{theorem}
Vectors $\eta_k,\Theta_k$ and $\Lambda_k$ denote the tensor components $\eta_{zk},\Theta_{zk}$ 
and $\Lambda_{zk}$. With the above definitions and results, the Frank and Burgers vectors are physical measures of defect which are given in terms 
of the sole strain and contortion tensors:
\begin{definition}
The Frank and Burgers vectors of surface $S$ are defined as
\begin{eqnarray}
\hspace{-85pt}\mbox{\scriptsize{MACROSCOPIC FRANK VECTOR}}\hspace{49pt}\Om_k(S)&:=&\int_S\Theta_k dS\label{frankmacro}\\
\hspace{-95pt}\mbox{\scriptsize{MACROSCOPIC BURGERS VECTOR}}\hspace{39pt}B_k(S)&:=&\int_S\Lambda_k  dS.\label{burgers1macro}
\end{eqnarray}
\end{definition}
As a consequence of (\ref{disclindens2}) and (\ref{frankmacro}) and Stokes theorem, the relation between the completed Frank tensor 
and the rotation gradient appears clear. Moreover, it results from (\ref{dislocdens2c}) and (\ref{burgers1macro}) that
$\left(\eth b\right)_{jk}:=\eth_j b_k$ appears instead of the displacement gradient. 

In the crystal dislocation-free regions 
(i.e. where the contortion vanishes), it results from the classical integral relation of infinitesimal elasticity 
that the multiple-valued rotation and displacement fields read
 \begin{eqnarray}
\omega_k&=&\omega_{0k}+\int_{x_0}^x\eth_m\omega_kd\xi_m
\label{rot}\\
u_i&=&u_{0i}-\epsilon_{ikl}\omega_k(x_l-x_{0l})
+\int_{x_0}^x\eth_mb_k(\xi)d\xi_m.\label{displ}
\end{eqnarray}

\subsection{Bravais metric and nonsymmetric connection as key geometrical objects}

A Riemannian metric is a smooth symmetric and positive definite tensor field $g_{ij}$. 
From its symmetry property, there is a smooth transformation $a_i^j$ such that 
$g_{ij}=a_i^ma_j^n\delta_{mn}$. The metric of the ``external observer'' on $\mathcal{R}(t)$ is the Euclidian metric $\delta_{ij}$. However, 
as soon as the macroscopic strain $\mathcal{E}_{ij}$ is given, another Riemannian metric can be defined on 
$\mathcal{R}(t)$, namely the
\begin{eqnarray}
\hspace{-84pt}\mbox{\scriptsize{BRAVAIS METRIC}}\hspace{115pt}g^B_{ij}=\delta_{ij}-2\mathcal{E}_{ij},\label{elastmetric1}
\end{eqnarray} 
where the term ``Bravais'' (from the notion of Bravais crystal \cite{KR90}) is to recall that it has not a purely elastic meaning. 

The use of this metric on defect-free regions of $\mathcal{R}(t)$ implies the 
existence of a one-to-one coordinate change between $\mathcal{R}(t)$ and $\mathcal{R}_0$, 
whose deformation gradient writes as $a_{mi}=g_{mn}^Ba^{n}_i=\delta_{mi}-\partial_i u_m$ where $u_m$ denotes a 
displacement-like field. Let us remark that since small displacements are considered, no distinction is to be made between upper and lower 
indices. 

In the presence of defects, the following object (which is said ``of anholonomity'' \cite{Schouten}) 
$\Omega_{ijk}:=\partial_ka_{ji}-
\partial_ia_{jk}$ is directly related to the strain incompatibility and hence does not vanish as soon as 
defects are present. This exactly signifies that there is no global system of coordinates $\{x^B_j(x_i)\}$ with a smooth 
transformation matrix $a_{ji}=\partial_i x^B_j$. In fact, such a smooth $a_{ij}$ -- or, 
equivalently, such a smooth displacement field only exist in the defect-free regions of the crystal. 

Quoting Cartan \cite{CART}, ``the Riemannian space is for us an 
ensemble of small pieces of Euclidian space, lying however to a certain degree amorphously'', while Kondo \cite{KONDO55} 
suggests that 
``the defective crystal is, by contrast [with respect to 
the above by him given
definition of perfect crystal], an aggregation of an immense number of small pieces of perfect crystals 
(i.e. small pieces of the defective crystal brought to their natural state in which the atoms are arranged on the regular positions of the perfect
 crystal) that cannot be connected with one other so as to form a finite lump of perfect crystals as an organic unity.''
\nl

From the elastic metric, we define the compatible symmetric Christoffel symbols 
\begin{eqnarray}
\hspace{-22pt}\mbox{\scriptsize{BRAVAIS CHRISTOFFEL SYMBOLS}}\hspace{35pt}\Gamma_{k;ij}^B=\frac{1}{2}\left(\partial_ig_{kj}^B+\partial_jg_{ki}^B-\partial_kg_{ij}^B\right),\label{compatconn}
\end{eqnarray} 
whose torsion $\Gamma_{k;[ij]}^B:=\Gamma_{k;ij}^B-\Gamma_{k;ji}^B$ vanishes, while its curvature
\begin{eqnarray}
\hspace{-25pt}\mbox{\scriptsize{BRAVAIS CURVATURE}}\hspace{50pt}R_{l;kmq}^B:=\left(\partial_q\Gamma_{l;km}^B+\tilde g^B_{np}\Gamma_{n;km}^B\Gamma_{p;lq}^B\right)_{[mq]}, \label{curv}
\end{eqnarray} 
with $\tilde g^B_{np}=\delta_{np}+\mathcal{E}_{np}$ the inverse of $g^B_{np}$ under the small strain assumption, and where symbol 
$[\cdot]$ denotes the skew symmetric index commutation operator (i.e., $A_{[mn]}=A_{mn}-A_{nm}$). In the terminology of Wang \cite{WANG67} 
and Noll 
 \cite{NOLL67} a connection such that $R_{l;kmq}^B$ vanishes is called \textit{a material connection}, while for 
such a connection $\Gamma_{k;[ij]}^B$ is denoted as \textit{the inhomogeneity tensor}.
\nl

Quoting Einstein, ``to take into account gravitation, we assume the existence of Riemannian metrics. 
But in nature we also have electromagnetic fields, which cannot be described by Riemannian metrics. 
The question arises: How can we add to our Riemannian spaces in a logically natural way an additional structure that 
provides all this with a uniform character ?'' 

In the present case, it is sufficient to replace gravitation by 
strain, and electromagnetic fields by line defects to paraphrase Einstein and raise the question of the apropriate  
connection inside the defective crystal. To be complete we should add that in order for the theory of dislocations 
to be closed, it should be combined with the theory of point defects which play a role at higher temperature (in the same way 
as Maxwell theory has to be combined with the theory of weak interactions, see Kr\"oner \cite{KR95}). 

The Bravais geodesics are those lines whose tangent vector $\tau_i$ is parallelly 
transported, hence solutions to $\tau_j\nabla^B_j\tau_i=0$, where $\nabla^B$ is the covariant derivative 
of $\Gamma^B$. It turns out that on these lines, the internal observer is not be able to recognize any defect line. 
Therefore, the above Bravais connection must be completed by a non-symmetric term.

The following geometric objects are introduced from the sole dislocation density 
(or equivalently by (\ref{dislocdens}) \& (\ref{dislocdens2c}) from the sole strain and contortion tensors):
\begin{definition}\label{defgeom}
\begin{eqnarray}
\hspace{-125pt}\mbox{\scriptsize{DISLOCATION TORSION:}}\hspace{100pt} T_{k;ij}&:=&-\frac{1}{2}\epsilon_{ijp}\Lambda_{pk}\label{torsion}\\
\hspace{-105pt}\mbox{\scriptsize{CONNECTION CONTORTION:}}\hspace{74pt}\Delta\Gamma_{k;ij}&:=&T_{j;ik}+T_{i;jk}-T_{k;ji}\label{contortion}\\
\hspace{-28pt}\mbox{\scriptsize{NON SYMMETRIC CHRISTOFFEL SYMBOLS:}}\hspace{23pt}\Gamma_{k;ij}&:=&\Gamma_{k;ij}^B-\Delta\Gamma_{k;ij}.
\label{connexion}
\end{eqnarray}
\end{definition}
According to Noll 
 \cite{NOLL67}  $\Delta\Gamma_{k;[ji]}$ is precisely the crystal \textit{inhomogeneity tensor} which will be shown in the following sections
to be directly related to the density of dislocations and disclinations.
\section{The macroscopic crystal as a non-Riemannian manifold}
By contrast with Kr\"oner's presentation, the present approach shows 
geometrical objects as defined from homogenization of mesoscopic  
measurable, objective physical fields (\ref{dislocdens2c}), (\ref{torsion}) \& (\ref{contortion}) 
whose identification with their physical macroscopic counterparts follows 
 as (proved) results. 

\subsection{Physical and geometrical torsions and contortions}
\label{sec:contortion}

The following lemma is easy to prove from the definitions.
\begin{lemma}\label{propgeom}
The tensor $g_{ij}^B$ defines a Riemannian metric. The symmetric Christoffel symbols $\Gamma_{k;ij}^B$ define a 
symmetric connection compatible with this metric, while $T_{k;ij}$ and $\Delta\Gamma_{k;ij}$ are skew-symmetric tensors 
w.r.t. $i$ and $j$ and $i$ and $k$, respectively. 
\end{lemma}
The following results makes the link between internal motion of the observer by parallel transport and the deformation and defect internal 
variables encountered.
\begin{theorem}[Physical and geometrical torsions]\label{gamma}
The Cristoffel symbols $\Gamma_{k;ij}$ define a nonsymmetric connection compatible with $g_{ij}^B$ whose 
torsion writes as $T_{k;ij}$.
\end{theorem}
\textbf{Proof.} 
It is easy to verify \cite{Dubrovin} that $\Gamma_{k;ij}$ is a connection since $\Gamma_{k;ij}^B$ is a 
connection and $\Delta\Gamma_{k;ij}$ is a tensor. Denoting by $\nabla_k$ (resp. $\nabla^B_k$) the covariant gradient 
w.r.t. $\Gamma_{k;ij}$ (resp. $\Gamma_{k;ij}^B$), and recalling that a connection is compatible with the metric 
$g_{ij}^B$ if the covariant gradient of $g_{ij}^B$ w.r.t. this connection vanishes, we find by (\ref{connexion}) that
\begin{eqnarray}
\nabla_{k}g_{ij}^B:&=&\partial_k g_{ij}^B-\Gamma_{l;ik}g_{lj}-\Gamma_{l;jk}g_{li}^B=\nabla^B_{k}g_{ij}^B+\Delta\Gamma_{l;ik}g_{lj}^B+\Delta\Gamma_{l;jk}g_{li}^B\label{tilde},
\end{eqnarray}
where in the RSH, the $1^{st}$ term vanishes by Lemma \ref{propgeom} while the $2^{nd}$ and $3^{rd}$ 
terms cancel each other since $\Delta\Gamma_{l;jk}g_{li}=\Delta\Gamma_{i;jk}=-\Delta\Gamma_{j;ik}$. It results 
that the connection torsion, i.e. the skew-symmetric part of $\Delta\Gamma_{j;ik}$ w.r.t. $i$ and $k$, writes as
\begin{eqnarray}
&&\hspace{-32pt}\frac{1}{2}\left(\Delta\Gamma_{j;ik}-\Delta\Gamma_{j;ki}\right)=-\frac{1}{2}\left(\Delta\Gamma_{i;jk}-\Delta\Gamma_{k;ji}\right)=\frac{1}{2}\bigl(\left(\Delta\Gamma_{k;ij}-\Delta\Gamma_{i;kj}\right)+\nonumber\\
&&\hspace{87pt}\left(\Delta\Gamma_{k;ji}-\Delta\Gamma_{k;ij}\right)-\left(\Delta\Gamma_{i;jk}-\Delta\Gamma_{i;kj}\right)
\bigr)\label{calcul}.
\end{eqnarray}
Observing that the $1^{st}$ term in the RHS side of (\ref{calcul}) writes as $\Delta\Gamma_{k;ij}$ while, by 
Definition \ref{defgeom} (Eq. (\ref{contortion})), the LHS and the two remaining terms of the RHS
of (\ref{calcul}) are equal to $T_{j;ik}, T_{k;ji}$ and $-T_{i;jk}$, respectively, the proof is complete.
{\hfill $\square$} 

\begin{theorem}[Physical and geometrical contortions]\label{deltagamma}
The connection contortion tensor $\Delta\Gamma_{k;ij}$ writes in terms of the dislocation contortion $\kappa_{ij}$ as
\begin{eqnarray}
\Delta\Gamma_{k;ij}=\delta_{k\kappa}\left(\delta_{i\alpha}\delta_{j\beta}\epsilon_{\kappa\alpha}\kappa_{z\beta}\right)
+\delta_{i\alpha}\delta_{jz}\epsilon_{\alpha\tau}\kappa_{\tau\kappa}&+&\delta_{iz}\delta_{j\beta}\epsilon_{\beta\tau}
\kappa_{\tau\kappa}-\delta_{kz}\delta_{i\alpha}\delta_{j\beta}\epsilon_{\alpha\beta}\kappa_{zz}.\nonumber
\end{eqnarray}
\end{theorem}
\textbf{Proof.} 
For $k=z$, by Definition \ref{defgeom}, the last statement of Lemma \ref{propgeom}, and (\ref{dislocdens2d}), 
it is found that $\Delta\Gamma_{z;ij}=\Delta\Gamma_{z;\alpha\beta}\delta_{i\alpha}\delta_{j\beta}$, 
with 
\begin{eqnarray}
\Delta\Gamma_{z;\alpha\beta}=T_{z;\alpha\beta}=-\frac{1}{2}\epsilon_{\alpha\beta}\Lambda_z
=-\frac{1}{2}\epsilon_{\alpha\tau}\delta_{\tau\beta}\Lambda_z
=\epsilon_{\alpha\tau}\kappa_{\tau\beta}.\nonumber
\end{eqnarray}
For $k=\kappa$, by Definition \ref{defgeom} and the last statement of Lemma \ref{propgeom}, it is found that
\begin{eqnarray}
\Delta\Gamma_{\kappa;ij}=\delta_{i\alpha}\delta_{j\beta}\left(T_{\kappa;\alpha\beta}+T_{\beta;\alpha\kappa}
+T_{\alpha;\beta\kappa}\right)+\delta_{i\alpha}\delta_{jz}T_{z;\alpha\kappa}+\delta_{iz}\delta_{j\beta}T_{z;\beta\kappa},
\nonumber
\end{eqnarray}
with $T_{z;\xi\kappa}=\epsilon_{\xi\tau}\kappa_{\tau\kappa}$
and $\displaystyle T_{\xi;\tau\nu}=-\frac{1}{2}\epsilon_{\tau\nu}\Lambda_\xi$. Since the combination of the terms 
in $\Theta_z$ vanish in $\Delta\Gamma_{\kappa;ij}$, the proof is completed by 
observing that $\epsilon_{\alpha\beta}\Lambda_\kappa+\epsilon_{\kappa\alpha}\Lambda_\beta=(\epsilon_{\alpha\kappa}
\epsilon_{\tau\nu})\epsilon_{\tau\beta}\Lambda_\nu=\epsilon_{\alpha\kappa}\Lambda_\beta=\epsilon_{\alpha\kappa}
\kappa_{z\beta}$.
{\hfill $\square$} 

In conclusion, the non-Riemannian crystal is described from a physical viewpoint by $\mathcal{E}$ 
and $\kappa$, that is by $15$ degrees of freedom. From a geometrical viewpoint the $15$ 
unknowns are the $6$ components of the symmetric Bravais metric and the $9$ (by (\ref{torsion}) \& (\ref{contortion}))
nonvanishing components of the connection contortion.

\subsection{The Bravais crystal}
\label{sec:bravais}

The following definition introduces two differential forms whose path integrations generalize (\ref{rot}) 
\&  (\ref{displ}) to the defective regions of the crystal.
\begin{definition}[Bravais forms]\label{formdiff}
\begin{eqnarray}
d\omega_j&:=\eth_\beta\omega_jdx_\beta,\label{domega}\\
d\beta_{kl}&:=-\Gamma_{l;k\beta}dx_\beta\label{dbeta}.
\end{eqnarray}
\end{definition}
In the literature the existence of an elastic macroscopic distortion field is generally postulated
together with the global distortion decomposition in elastic and plastic parts (for a rigorous justification 
of the latter, see \cite{KR96}). 
The present approach renders however possible to avoid this a-priori decomposition. Nevertheless, the following theorem introduces 
rotation and distortion fields (which must not be identified with the rotation and distortion as related to the macroscopic strain) 
in the absence of disclinations. As a consequence and in contrast with the classical literature where it is basically postulated that 
dislocation density is the distortion curl, this relationship is here proved.
\begin{theorem}\label{elasticrotdist}
If the macroscopic disclination density vanishes, there exists rotation and distortion fields defined as
\begin{align}
&\hspace{-50pt}\mbox{\scriptsize{BRAVAIS ROTATION}}\ &\omega_j(x)&:=\omega_j^0+\int_{x_0}^x d\omega_j,\label{omega}\\
&\hspace{-50pt}\mbox{\scriptsize{BRAVAIS DISTORTION}}\ &\beta_{kl}(x)&:=\mathcal{E}_{kl}(x^0)-\epsilon_{klj}\omega^0_j
+\int_{x_0}^x d\beta_{kl},\label{beta}\\
&&&=\mathcal{E}_{kl}(x)-\epsilon_{klj}\omega_j(x)
\end{align}
where $\omega_j^0$ is arbitrary and the integration is made on any line with endpoints $x_0$ and $x$. Moreover,
\begin{eqnarray}
\partial_\alpha\beta_{k\beta}=\partial_\alpha\mathcal{E}_{k\beta}+\epsilon_{kp\beta}\eth_\alpha\omega_p\quad\mbox{and}
\quad
\epsilon_{\alpha\beta}\partial_\alpha\beta_{k\beta}=\Lambda_{zk}.\label{decompmacro}
\end{eqnarray}
\end{theorem}
\textbf{Proof.} 
By Definition \ref{defgeom}, the symmetric part of the connection writes as
\begin{eqnarray}
-\Gamma_{(l;k)\beta}dx_\beta=-\frac{1}{2}\partial_\beta g_{kl}^B dx_\beta=-\frac{1}{2}\partial_m g_{kl}^B dx_m
=\partial_m \mathcal{E}_{kl}dx_m=d\mathcal{E}_{kl},\nonumber
\end{eqnarray}
while, by Definition \ref{defgeom} and Theorem \ref{deltagamma}, the skew-symmetric part writes as
\begin{eqnarray}
-\Gamma_{[l;k]\beta}&=&-\frac{1}{2}(\partial_k g_{l\beta}^B -\partial_l g_{k\beta}^B )+\Delta\Gamma_{l;k\beta}=
\partial_k\mathcal{E}_{l\beta}-\partial_l\mathcal{E}_{k\beta}+\Delta\Gamma_{l;k\beta}.\nonumber
\end{eqnarray}
Observing, by (\ref{disclindens}) and Definition \ref{formdiff} and Theorem \ref{deltagamma}, that
$\displaystyle d\omega_j=\eth_\beta\omega_j dx_\beta$ $=-\frac{1}{2}\epsilon_{lkj}\Gamma_{[l;k]\beta}dx_\beta$,
it results that $d\beta_{kl}=d\mathcal{E}_{kl}-\epsilon_{klj}d\omega_j$.
Under the assumption of a vanishing macroscopic disclination density, 
the existence of single-valued Bravais rotation and distortion fields follows from (\ref{domega}), 
(\ref{disclindens2}) \& (\ref{frankmacro}). Moreover, since 
$\partial_\alpha\beta_{k\beta}=\partial_\alpha\mathcal{E}_{k\beta}-\epsilon_{k\beta j}\eth_\alpha\omega_j$, 
Eq. (\ref{decompmacro}) is satisfied by (\ref{disclindens2}) \& (\ref{dislocdens2c}).
{\hfill $\square$} 
\begin{remark}
Eq. (\ref{omega}) indicates that symbol $\eth$ in (\ref{domega}) becomes a true derivation operator in the absence of disclinations. 
\end{remark}

\begin{remark}\label{seulconn}
Referring to ``Bravais'' instead of ``elastic'' rotation and distortion fields is devoted to highlight that these 
quantities do not have a purely elastic meaning. In fact, the Bravais metric is not even needed since the internal observer
only requires the prescription of the connection, and subsequent path integration of the forms:
\begin{eqnarray}
d\mathcal{E}_{kl}=-\Gamma_{(l;k)\beta}dx_\beta,\quad d\omega_j:=-\frac{1}{2}\epsilon_{lkj}\Gamma_{[l;k]\beta}dx_\beta\quad\mbox{and}\quad
d\beta_{kl}:=-\Gamma_{l;k\beta}dx_\beta\nonumber.
\end{eqnarray}
\end{remark}

\begin{remark}\label{remark5}
The Bravais distortion does not derive from any ``Bravais displacement'' in the presence of dislocations. 
In fact, around a closed loop $C$, even if the disclination density vanishes, the  differential of the displacementas 
$du_k:=\beta_{k\alpha}dx_\alpha$ verifies by Theorem \ref{elasticrotdist}:
\begin{eqnarray}
\int_C du_k=\int_S\epsilon_{\alpha\beta}\partial_\alpha\beta_{k\beta}dS
=\int_C \eth_\beta b_k dx_\beta=\Lambda_k(S)=\int_S\epsilon_{\alpha\beta}\partial_\alpha\eth_\beta b_k dS.
\end{eqnarray}
\end{remark}

\begin{remark}\label{remark4}
Theorem \ref{gamma} defines an operation of parallel displacement according to the Bravais lattice geometry. 
The parallel displacement of any vector $v_i$ along a curve of tangent vector $dx^{(1)}_\alpha$ is such that 
$dx^{(1)}_\alpha\nabla_\alpha v_i=0$ and hence that the components of $v_i$ vary according to the 
law $d^{(1)}v_i=-\Gamma_{i;j\beta}v_j dx^{(1)}_\beta$ \cite{Dubrovin}. This shows the macroscopic Burgers 
vector and dislocation density together with the Bravais rotation and distortion fields as reminiscences of the 
defective crystal properties at the atomic, mesoscopic and continuum scales. In fact, if $dx^{(1)}_\nu, dx^{(2)}_\xi$ are two infinitesimal vectors 
with the associated area $dS:=\epsilon_{\nu\xi}dx^{(1)}_\nu dx^{(2)}_\xi$, it results from Eq. (\ref{torsion}) and the skew symmetry of $\Gamma_{k;\alpha\beta}$ that, in the absence of disclinations, 
\begin{eqnarray}
dB_k=\Lambda_{z\alpha} dS
=-\epsilon_{\alpha\beta}\Gamma_{k;\beta\alpha}dS
=-\Gamma_{k;\beta\alpha}(dx^{(1)}_\alpha dx^{(2)}_\beta-dx^{(1)}_\beta dx^{(2)}_\alpha)\nonumber,
\end{eqnarray}
whose right-hand side appears as a commutator verifying the relation
\begin{eqnarray}
dB_k=\epsilon_{\alpha\beta}\partial_\alpha\eth_\beta b_k dS=-\epsilon_{\alpha\beta}d^{(\alpha)}(dx^{(\beta)})\nonumber.
\end{eqnarray}
\end{remark}

\subsection{Motion of the internal observer}
\label{sec:int}

The internal observer will be represented by the $\textbf{k}^{th}$ geodesic  basis element ${\bf{e_k}}(x)$ solution to
\begin{eqnarray}
({\bf{e_k}})_j\nabla_j({\bf{e_k}})_l=0\quad\mbox{(with no summation on $\texttt{\textbf{k}}$)}, 
\end{eqnarray}
where $\nabla$ is the covariant derivative of $\Gamma$ as given by (\ref{connexion}). We have seen 
in the above two sections that it was sufficient to provide him with a connection, i.e. with a law of parallel transport inside
the crystal. In fact at this stage the internal observer is not able to measure distances, while he can measure the disclination
 (resp. dislocation) content of a 
surface $S$ by boundary measurements of $\eth_\beta b_k$ (resp. $\eth_\beta\omega_k$) on the curve $C$ enclosing $S$ 
(which depend merely on $\Gamma$ -- cf Remarks \ref{seulconn}-\ref{remark4}).

The notion of metric connection can be explained as follows. Let the external observer be equipped with the 
Bravais metric and Cartesian coordinate system $\{x_i\}$. Since $\Gamma_{l;km}=\nabla_m ({\bf{e_k}})_l$, we have 
on a portion $A-B$ of geodesic $\texttt{k}$,
\begin{eqnarray}
({\bf{e_k}})_l(B)-({\bf{e_k}})_l(A)=\int_A^B\Gamma_{l;km}dx_m
=\lim_{N\to\infty}\sum_{1\leq i\leq N}\Gamma_{l;km}(x^i)({\bf{e_k}})_m(x^i)\Delta s^i,\nonumber
\end{eqnarray} 
where $x^i$ are discretisation points on the curve with endpoints $x^1=A$ and $x^N=B$, and $\Delta s^i$ a tending to zero element of the geodesic. 
Moreover, if the connection is compatible with the metric 
$g^B$, the angles between these lattice vectors and their (unit) length remain invariant during 
parallel transport. So, we understand Kr\"{o}ner \cite{KR92} when he says: ``when a lattice vector is parallelly displaced using 
$\Gamma$ along itself, say $1000$ times, then its start [say, $A$] and goal [say, $B$] are separated by 1000 atomic spacings, as 
measured by $g^B$. Because the result of the measurement by parallel displacement and by counting lattice steps 
is the same, we say that the space is metric with respect to the connection $\Gamma$.'' 

Moreover, as long as  $({\bf{e_k}})_l(A)$ (the internal observer) is parallely transported along a closed curve $C$ with start- and endpoint 
$A$, the gap as created when he comes back to his orgin can be measured by the external observer, since by Stokes theorem
\begin{eqnarray}
&&({\bf{e_k}}^\shortparallel)_l-({\bf{e_k}})_l=\int_C\Gamma_{l;km}dx_m=\int_{S}\epsilon_{pqm}\partial_q
\Gamma_{l;km}dS_p\nonumber\\
&=&\int_{S}\epsilon_{pqm}\left(\nabla_q
\Gamma_{l;km}+\left(\Gamma_{l;pm}\Gamma_{p;kq}+\Gamma_{l;pq}\Gamma_{p;km}\right)+\Gamma_{l;kp}\Gamma_{p;mq}\right)dS_p
\nonumber,
\end{eqnarray} 
where $({\bf{e_k}}^\shortparallel)_l$ denotes the  base $({\bf{e_l}})_k$ after being parallely transported along $C$. Since the term inside the 
parenthesis is symmetric in $m$ and $q$, we have \cite{Dubrovin}:
\begin{eqnarray}
=\int_{S}\epsilon_{pqm}\frac{1}{2}\left(\nabla_{[q}
\nabla_{m]} ({\bf{e_k}})_l+\nabla_p({\bf{e_k}})_l T_{p;mq}\right)dS_p=\int_{S}\epsilon_{pqm}\frac{1}{2}R_{l;nmq}({\bf{e_k}})_ndS_p,\label{internmotion}
\end{eqnarray} 
with the definition of the Riemannian curvature tensor
\begin{eqnarray}
R_{l;kmq}:=R_{l;kmq}^B+\Delta R_{l;kmq}, \label{curv2}
\end{eqnarray}
where by (\ref{curv}) \& (\ref{connexion}), $R^B$ and $\Delta R$ denote the Riemann curvature tensors 
associated to $\Gamma^B$ and 
$\Delta\Gamma$, respectively. 
\nl

By (\ref{internmotion}), the internal observer is convinced to return to his startpoint
while the external observer however can see the gap as created by the crystal curvature, itself resulting from the presence of defects.

\subsection{Geometric and physical curvatures}
\label{sec:curv}

Let us remark that in the absence dislocations ($T=\Delta R=\Lambda=0$), the gap is merely due to curvature 
effects with a curvature tensor directly related by (\ref{etakmacro}) to the disclination density by 
$R_{l;kmq}=-\epsilon_{lki}\epsilon_{mqj}\Theta_{ij}$. It should however be noted that in the absence of disclinations, the curvature 
is not vanishing but depends on the sole contortion, since from (\ref{etakmacro}) \& (\ref{curv2}),
\begin{eqnarray}
R_{l;kmq}=-\epsilon_{lki}\epsilon_{mqj}\epsilon_{ipn}\partial_p\kappa_{jn}+\Delta R_{l;kmq},\label{curv4}
\end{eqnarray} 
where by Theorem \ref{deltagamma}, $\Delta R$ is linearly related to the contortion. It is computed from (\ref{curv4}) that the Ricci and Gauss 
curvatures \cite{Dubrovin} read
\begin{eqnarray}
 \hspace{-68pt}\mbox{\scriptsize{RICCI CURVATURE}} \hspace{38pt}R_{kq}^B&:=&R_{l;kmq}^B=R_{p;kpq}^B=\eta_{kq}-\delta_{kq}\eta_{pp}\label{Ricci}\\
 \hspace{-65pt}\mbox{\scriptsize{GAUSS CURVATURE}} \hspace{36pt}R^B&:=&\frac{1}{2}R^B_{pp}=-\eta_{pp},\label{Gauss}
\end{eqnarray}
while Einstein tensor reads
\begin{eqnarray}
-\frac{1}{4}\epsilon_{lki}\epsilon_{mqj}R_{l;kmq}^B=\eta_{ij}=R_{ij}^B-\delta_{ij}R^B
\end{eqnarray}
in the presence of dislocations and disclinations, thereby contradicting Kr\"oner who identified Einstein 
tensor with the disclination density in \cite{KR92}.
 
Moreover, since the macroscopic strain can be decomposed into (symmetric) compatible and (symmetric) solenoidal 
parts \cite{VGD2009}, where only the second one as denoted by $\mathcal{E}^s$ is relevant for the incompatibility tensor, it results that 
its trace 
$\mathcal{E}_{pp}^s$ satisfies by (\ref{Gauss}) $-\Delta \mathcal{E}_{pp}^s=R^B$, thereby
showing how the Gauss curvature is related to the variation of matter density. 

\subsection{Summary of the non-Riemannian metric crystal}
\label{sec:sum}
The crystal equiped with
$\{g^B,\Gamma\}$ has the following properties: (i) the geodesics of $\Gamma$ are the crystallographic lines; 
(ii) the effect of parallel displacement 
of the internal observer (equiped with $\Gamma$) along a crystallographic line is equivalent to counting the lattice steps;  
(iii) the defect content, i.e. disclination and/or dislocation densities can be computed from 
measures of $\Gamma$ only; (iv) the torsion of $\Gamma$ is merely due to the presence of dislocations, while its curvature is 
due to the presence of both disclinations and dislocations; 
(v) in the absence of disclinations, there exists a single-valued rotation and distortion field; 
(vi) if and only if there are no defect lines, $\Gamma$ is Euclidean and there exists a holonomic 
coordinate system. In the latter case only, one can properly speak of a reference configuration, 
of single-valued rotation, displacement and distortion fields, with the macroscopic strain compatible with the displacement field.
\nl

\begin{figure}[htbp]
\begin{center}
\includegraphics[width=11cm,height=7cm]{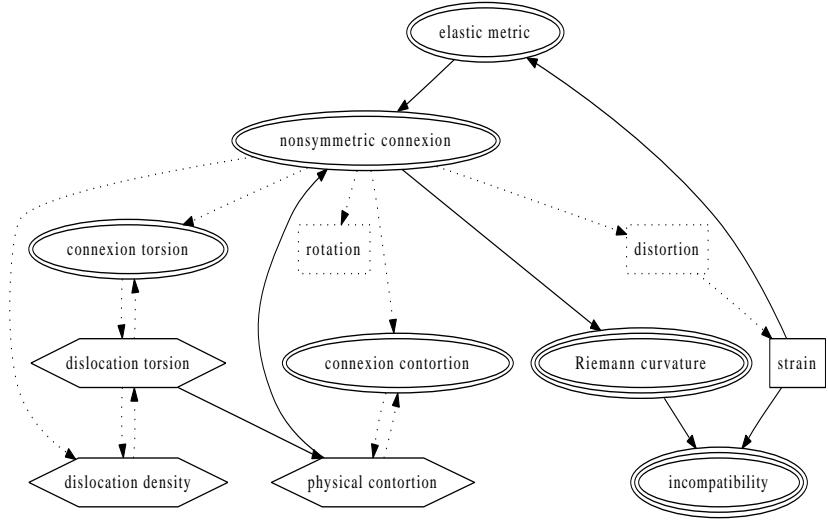}
\end{center}
\caption{Link between physics and geometry of defects}\label{geom}
\end{figure}

Figure \ref{geom} illustrate the inseparable link between physics and geometry. On the one hand, the physical fields can be set apart: the deformation 
and defect internal variables are shown in rectangular and hexagonal boxes, respectively. On the other hand, the purely geometrical object are 
in oval boxes. The double boundary line means that the quantity contains differential combinations of other fields (as connected by arrows), 
while single lines mean algebraic combinations only.

The main deformation field is the strain, while the distortion and rotation are only obtained in the absence of disclinations (see Theorem 
\ref{elasticrotdist}). Because they depend on an arbitrary point where there value is assumed known, they are considered inappropriate  
as model variables. Concerning defect 
internal variables, one could indifferently chose the dislocation torsion or contortion. Let us mention that since strain instead of distortion is
chosen, the deformation and defect variables should be considered as independent physical fields.

\subsection{Nonmetricity, teleparallelism and the paradox of the flat crystal}
\label{sec:distpar}
Let us first remark that the notions of metric and of connection must be considered as distinct. This has been
emphasized in \cite{BOURG} where it is recalled that historically it has not been so for a long time (including Einstein literature).
Here, we have seen that the metric is attached to the notion of external observer, while the connection is attached to the notion of
parallel displacement of the internal observer inside the crystal. We have seen that parallel displacement with $\Gamma^B$ and 
counting-step measurements are the same in a crystal filled with line defects. This was true because $\Gamma^B$ was compatible with 
$g^B$ and hence crystallographic basis elements remain crystallographic as transported along the crystallographic lines. 

Suppose now that the crystal also contains point defects. According to Kr\"oner \cite{KR96}, ``nonmetricity means that length measurements are disturbed. It is easy to see that this just occurs in the presence of 
point defects. In fact, when counting atomic steps along crystallographic lines to measure distance between two atoms, [the internal observer] feels 
disturbed when suddently a vacancy or an intersticial emerges instead of another atoms [of the perfect crystal].''

Let $C_V$ and $C_I$ be the scalar vacancy (resp. interstitial) concentration, that is the number of vacancies (resp. interstitials) per unit 
volume of crystal. Then, the following metric as proposed by Kr\"oner \cite{KR90}:
\begin{equation}
 g'=(1+C_I-C_V)^2g^B,\label{newmetric}
\end{equation} 
verifies $dV=\sqrt{\det g'} dV_0=(1+\Delta C)\sqrt{\det g^B} dV_0$, with $dV_0$ the volume element of the stress-free crystal and $dV$ that of 
the actual one, and where $\Delta C=C_I-C_V$ is the excess atomic content of $dV$. An evolution equation for $C_V$ and $C_I$ (and hence for $\Delta C$ and $g'$) 
will be given in Section \ref{concl}.

It is clear that the non-metricity defined as $Q_{j;ik}:=\hat \nabla_j g'_{ik}\neq 0$ \cite{Schouten,KR92} must now enter the geometric point- and 
line-defect model. Differentiation $\hat\nabla$ is here intended with respect to connection $\hat\Gamma$, as defined by
\begin{eqnarray}
\mbox{\scriptsize{CHRISTOFFEL SYMBOLS WITH POINT DEFECTS:}}\hspace{4pt}\hat\Gamma_{k;ij}:=\Gamma'_{k;ij}-\Delta\Gamma_{k;ij}
-\frac{1}{2}\delta\Gamma_{k;ij}\label{connexion2}
\end{eqnarray}
with $\Delta\Gamma_{k;ij}$ given by (\ref{torsion}) \& (\ref{contortion}) and where $\Gamma'_{k;ij}:=\frac{1}{2}\left(\partial_ig'_{kj}+
\partial_jg'_{ki}-\partial_kg'_{ij}\right)$ and
\begin{eqnarray}
\mbox{\scriptsize{NONMETRIC CONTORTION:}}\hspace{100pt}\delta\Gamma_{k;ij}:=Q_{j;ik}+Q_{i;jk}-Q_{k;ji}.\label{contortion2}
\end{eqnarray}
Since $Q$ is a tensor quantity, it is expected to play a role for the physical description of the crystal (and to obey an evolution equation as related 
to the other defects). 
\nl

The paradox of the flat 
crystal is the fact that a defective solids with in addition to defect lines a certain amount of point defects
can recover a vanishing curvature if the following balance holds:
\begin{eqnarray}
\hspace{-35pt}\mbox{\scriptsize{TELEPARALLELISM ASSUMPTION:}}\hspace{30pt}\delta R_{l;kmq}=-(R'_{l;kmq}+\Delta R_{l;kmq}), \label{curv6}
\end{eqnarray}
where the three terms are the curvature of $\delta\Gamma,\Gamma'$ and $\Delta\Gamma$, respectively. This is what Kr\"oner 
(and Bilby et al \cite{BILBY}) calls teleparallelism, by 
this meaning that the global connection curvature vanishes and hence that the internal observer ends up parallel when travelling along a loop. It 
should be
emphasized that teleparallelism is often considered as a working assumption \cite{BILBY,KR96}. 

However, we rather follow Kr\"oner \cite{KR90} when he says 
that 
``curved crystals are possible only if the curvature is, in some sense, compatible with the considered crystal structure'', which means that instead 
of a 
flat crystal given by (\ref{curv6}), the connection curvature $\hat R$ of the actual crystal should be such that point and line defects accomodate to 
satisfy
 \begin{eqnarray}
\hat R_{l;kmq}=\delta R_{l;kmq}+(R'_{l;kmq}+\Delta R_{l;kmq}), \label{curv7}
\end{eqnarray}
Particularizing (\ref{curv7}) we learn from identity $\hat R_{(l;k)mq}+\hat\nabla_{[m}Q_{q];lk}+T_{p:mq}Q_{p;lk}=0$ \cite{Schouten,KR92} that point defects 
and dislocations must be geometrically related, which phenomena is well known from solid-state physicist: ``dislocations moving perpendicular to their Burgers vector produce point defects, and similar processes occur when 
dislocations cut each other'' \cite{KR90} (concerning non metricity, see also \cite{BENABR}).

\section{Concluding remarks: the choice of model variables}\label{concl}

Let us conclude by attempting to answer Kr\"oner's question of the Introduction. To the knowledge of the author, this question has not 
been answered yet, or to say the least, there is still no agreement on the answer. 

In fact it depends of the physics which one wants to capture. If motion of dislocations is modelled, then not
only the conservative glide but also the non-conservative climb mode must be taken into account. Non conservation is due to the presence, 
creation, anihilation and motion of point defects, and these processes require high temperatures and non-negligible temperature gradients 
\cite{VGetal,Phil}. Therefore a complete model of dislocation motion must be thermodynamical, away from thermal equilibrium (i.e. irreversible), 
and coupled with the motion of point defects. In \cite{VGetal} the following set of 
PDEs showed very good results to model point defects:
\begin{eqnarray}
\frac{DC_K}{Dt}&=&\nabla\cdot\left(D_K\nabla C_K+\tilde D_KC_K\nabla T\right)-P,\label{eqci}
\end{eqnarray}
with the Lagrangian derivative $D/Dt$, and where $C_K, D_K, \tilde D_K$ and $P$ mean (scalar) concentration, (tensorial) equilibrium diffusion, 
thermodiffusion 
and (scalar) recombination ($K=I$ for interstitials and $K=V$ for vacancies).

Concerning the motion of dislocations (we assume here that disclinations are negligible), a similar equation as (\ref{eqci}) should be proposed, 
with an inter-dislocation recombination term and a term of interaction with point defects, collectively denoted by $\tilde P$ (also appearing in (\ref{eqci})). However, 
the 
dislocation density cannot be scalar (which is the case in most of the current models available in the literature) but must be the tensorial 
$\Lambda$ 
(or equivalently the contortion $\kappa$). The PDE could read
\begin{eqnarray}
\frac{D\kappa}{Dt}&=&\nabla\cdot\left(D\nabla \kappa+\tilde D\kappa\nabla T\right)-\tilde P,\label{eqdisloc}
\end{eqnarray}
with appropriate boundary condition and where $D$ and $\tilde D$ are tensor diffusivities of order $4$. Moreover the contortion verifies the 
conservation law \cite{VGD2009} 
$\nabla\cdot\kappa=\nabla\left(\texttt{tr}\ \kappa\right)$, meaning that the mesoscopic dislocations are loops or end at the crystal boundary, in such a way that Eq. (\ref{eqdisloc}) amounts to a system of $6$ coupled PDEs. 
Moreover the expression of $\tilde P$ must somehow satisfy the geometric interaction between defects as given by (\ref{curv7}). 
\nl

\begin{figure}[htbp]
\begin{center}
\includegraphics[width=11cm,height=6.5cm]{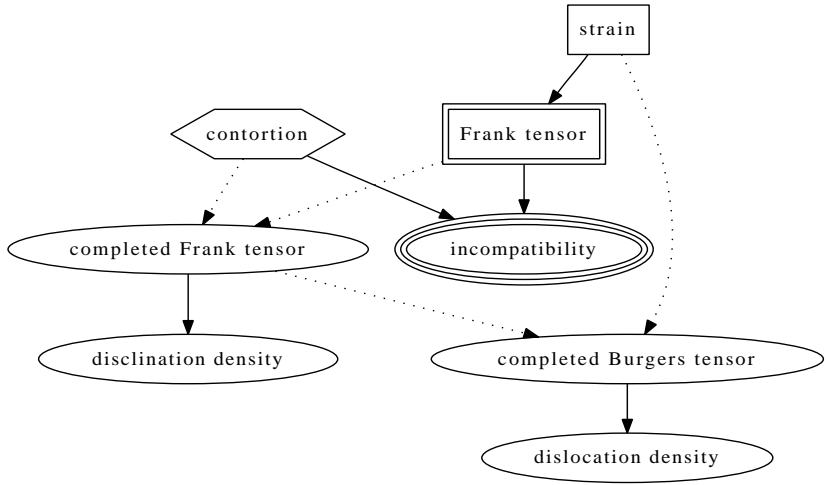}
\end{center}
\caption{The deformation and defect state variables}\label{physics}
\end{figure}

Concerning the deformation variables, let us first observe Figure \ref{physics}. 
Incompatibility is the final quantity as obtained by recursive
differentiation of either the strain (twice), or the contortion (once). It hence shows the ultimate convergence of the initially set apart deformation 
and defect variables. The other oval boxes denote defect variables obtained from the two key model variables: strain and contortion. 

It should be
observed that all relations between strain $\mathcal{E}$, Frank tensor $\overline\partial\omega$, contortion $\kappa$ and incompatibility are obtained by means of recursive application of the curl 
operator (either $\nabla\times$ or $\times\nabla$). 
\nl

As a first step, Kr\"oner 
proposed an (``athermal'') Gibbs free energy reading $W=\tilde W(F^e,\Lambda)$ with $F^e$ the elastic deformation gradient. 
Restricting hence to statics, he thereby attempted to answer his question \cite{KR95}: ``what are the independent (extensive) [explicit state] 
variables entering the free energy (at constant temperature)?'' However in \cite{KR96}, he recognizes that the 
use of $F^e$ is inevitably ambiguous because the elasto-plastic decomposition is not unique. 

According to our theory, the free energy 
naturally reads from the diagram on Fig. \ref{physics} as a first strain gradient (see \cite{LEONEP,LAZMAUG}) model $W=W_1(\mathcal{E},
\overline\partial\omega;\kappa)$ where the strain gradient is however replaced by its curl. 
Equivalently it could read $W=W_2(\mathcal{E},\eth\omega,\Lambda)$ by combination of the last two
variables of $W_1$, or even 
$W=W_3(\mathcal{E},\eth b,\Lambda)$ by combination of the first two variables of $W_2$. Let us observe that according to 
Theorem \ref{elasticrotdist}, $W_3$ can nevertheless be compared with $\tilde W$ as soon as $\eth b$ is identified with a distortion (i.e. a deformation gradient), 
although not becessarily the elastic one (cf. Remark \ref{seulconn}). Moreover,  a curl differential relation is also observed between $W_3$ last two 
variables (cf. Remark \ref{remark5}). 
A recent thermodynamic analysis with $W_3$ has been remarkably reported by Berdichevsky \cite{BERD08} where 
$\eth b$ is identified with the plastic distortion. Let us remark however that by (\ref{disclindens}) 
\& (\ref{dislocdens}), $\eth b$ is not, because of the prescription of the arbirary $x_0$,
an unambiguous state variable, as opposed to $\eth \omega$.

They are however reasons to be tempted by the choice $W=W_1$ because (i) all variables are explicit state variables defined by objective fields 
(which do not appeal to reference configurations, arbitrary plastic parts, or points $x_0$), (ii) there is a distinction between (strain-like) deformation 
 and (internal) defect variables, (iii) all variables have clear and unambiguous physical and a geometrical meanings. If needed, all other variables (such as $\nabla\mathcal{E},
\eth\omega,\eth b,\Lambda,\eta$) can be recovered 
as implicit state variables of the model \cite{KR95}. Moreover, if applying the curl operator twice to the strain and once to the contortion, then 
the only additional model variable  naturally appearing is the incompatibility, of both deformation and defect nature. 
This sounds like a closure on the recursive iteration for (higher-order) models.

Nonetheless, we rather prefer to introduce incompatibility through Kr\"oner's formula (\ref{etakmacro}) as a constraint to the Gibbs energy 
with $6$ degrees of freedom as coupling strain, Frank tensor and contortion:
\begin{eqnarray}
 \nabla\times\mathcal{E}\times\nabla=\nabla\times\overline\partial\omega=\kappa\times\nabla, \label{KR}
\end{eqnarray}
while the equilibrium law writes as the following equation with $3$ d.o.f.:
\begin{eqnarray}
 -\nabla\cdot\sigma=f\quad\mbox{with}\quad  \sigma:=\frac{\partial W}{\partial\mathcal{E}}\quad\mbox{and}\quad W=\mathcal{W}(T,\nabla T; \mathcal{E},\overline\partial\omega;\kappa),\label{equ}
\end{eqnarray}
where $f$ is the sum of external forces and of configurational (internal) pseudo-forces directly related to $\kappa$ and to the
derivative of the so-called dislocation moment stress $\partial W_1/\partial\kappa$ \cite{KR96,BERD08}. Let us also mention that the 
additional constraint of incompressibility must be added in order to avoid climb and point defects \cite{KR95}. Also, a remarkable discussion on the 
nature of $W=W_4(\eth b)$ (with identification of $\eth b$ with the distortion) in a nonlinear and variational setting can be found 
in \cite{PAL2008}.
\nl

Summarizing, an athermal model of dislocations requires to solve equations (\ref{eqdisloc})-(\ref{equ}) which involve a total of $15$ degrees of feedom. This number
is exactly the number of d.o.f. required by the internal observer to parallel displace inside the crystal (through the nonsymmetric connexion). 
Moreover to be closed the theory must involve point and line defects, and hence must consider high temperature and temperature gradients. 
So, as recognized by Kr\"oner \cite{KR96}, the dislocation model must be cast within the general frame of irreversible thermodynamics because the 
time variations of the internal variables create thermal dissipation.

This is the main reason why a huge work remains to be done in order to determine, e.g., the stress-strain relation, all other constitutive laws, 
the diffusion coefficients (which depend on the crystal internal symmetries, glide planes, etc), and the defect interaction/production terms.
\nl

The author is unable to answer definitely any of these questions but will pursue research on the topic. 
This paper aims at recalling Ekkehart Kr\"oner's 
legacy, and in particular the fundamental questions he raised which are still open and crucial nowadays.  It is also aimed at stressing that solutions 
to dislocation modelling will most probably arise from a strong interplay between mathematics and physics, as remarkably done by Kr\"oner along his papers 
\cite{KR55}-\cite{KR01}.


\begin{thebibliography}{99}
\bibitem{AFP2000}
     \newblock L. Ambrosio, N. Fusco and D. Palara, 
     \newblock ``Functions of bounded variation and free discontinuity problems," 
     \newblock Oxford Mathematical Monographs,
     \newblock Oxford, 2000.

\bibitem{ANT70a}
     \newblock K. H. Anthony,
     \newblock \emph{Die Reduktion von nichteuklidischen geometrischen Objekten in eine euklidische Form und physikalische 
Deutung der Reduktion durch Eigenspannungszust\"ande in Kristallen},
     \newblock Arch. Rational Mech. Anal., \textbf{37}, 3,  (1970), 43--88.

\bibitem{ANT70b}
     \newblock K. H. Anthony,
     \newblock \emph{Die Theorie der Disklinationen},
     \newblock Arch. Rational Mech. Anal., \textbf{39}, 1,  (1970), 161--180.

\bibitem{BENABR} 
     \newblock S. Ben-Abraham,
     \newblock \emph{Generalized stress and non-Riemannian geometry},
    \newblock in ``Fundam. Aspects of Dislocation Theory" (Nat. Bur. Stand. (U.S.)),
                Spec. Publ \textbf{317}, II, (1970), 943--962.

\bibitem{BERD08}
     \newblock V. L. Berdichevsky,
     \newblock \emph{Continuum theory of dislocations revisited},
    \newblock Continuum Mech. Thermodyn., \textbf{18},  (2006), 195--222.

\bibitem{BILBY}
     \newblock B. A. Bilby, R. Bullough and E. Smith,
     \newblock \emph{Continuous distribution of dislocations: a new application of the methods of non-Riemannian geometry},
    \newblock Proc. Roy. Soc. London A, \textbf{231}, 1,  (1955), 263--273.

\bibitem{BOURG}
    \newblock J.-P. Bourguignon,
    \newblock  \emph{Transport parall\`{e}le et connexions en G\'{e}om\'{e}trie et en Physique. (French)[Parallel transport and connections in geometry and physics]  
1830--1930: a century of geometry (Paris, 1989)},
    \newblock  in Lecture Notes in Phys. \textbf{402}
    \newblock       Springer, Berlin,   (1992), 150--164.

\bibitem{BUCAEP}
    \newblock I. Bucataru and M. Epstein
    \newblock \emph{Geometrical theory of dislocations in bodies with microstructure},
    \newblock Journal of Geometry and Physics \textbf{52}, (2004), 57--73.

\bibitem{CART} 
     \newblock E. Cartan,
     \newblock \emph{Sur une generalisation de la notion de courbure de Riemann et les espaces a torsion},
     \newblock C. R. Acad. Sci. Paris, \textbf{174} (1922), 593--597.

\bibitem{LEONEP}
     \newblock M. de Le\'{o}n and M. Epstein,
     \newblock \emph{The geometry of uniformity in second-grade elasticity},
     \newblock Acta Mechanica \textbf{114}, (1996), 217--224.

\bibitem{LEONEP2}
     \newblock M. de Le\'{o}n and M. Epstein,
     \newblock \emph{Geometrical Theory of Uniform Cosserat Media},
     \newblock Journal of Geometry and Physics \textbf{26}, (1998), 127--170.

\bibitem{Dubrovin}
     \newblock B. A. Dubrovin, A. T. Fomenko, S. P. Novikov, 
     \newblock ``Modern geometry - methods and applications," 
     \newblock 2nd edn.,
     \newblock Springer-Verlag, New York, 1992.

\bibitem{EP2010}
     \newblock M. Epstein,
     \newblock ``The geometrical language of continuum mechanics," 
     \newblock Cambridge University Press, Cambridge, 2010.

\bibitem{EPMAUG}
     \newblock M. Epstein and G. A. Maugin, 
     \newblock \emph{The energy-momentum tensor and material uniformity in finite elasticity},
     \newblock Acta Mech. \textbf{83}, (1990), 127--133.
 
\bibitem{KONDO52}
    \newblock K. Kondo,
    \newblock  \emph{On the geometrical and physical foundations of the theory of yielding},
    \newblock  in ``Proc. 2nd Japan Nat. Congr. Applied Mechanics,"
               Tokyo  (1952), 41--47.

\bibitem{KONDO55}
    \newblock K. Kondo,
    \newblock  \emph{Non-Riemannian geometry of the imperfect crystal from a macroscopic viewpoint},
    \newblock  in ``RAAG Memoirs of the unifying study of basic problems in engineering sciences by means of geometry,"
                Vol.1, Division D, Gakuyusty Bunken Fukin-Day, Tokyo  (1955), 458--469.

\bibitem{KR55}
     \newblock E. Kr\"oner,
     \newblock \emph{Die Spannungsfunktionen der dreidimensionalen anisotropen Elastizit\"atstheorie},
     \newblock Z. Physik, \textbf{140} (1955), 386--398.

\bibitem{KR60}
     \newblock E. Kr\"oner,
     \newblock \emph{Allgemeine Kontinuumstheorie der Versetzungen und Eigenspannungen},
     \newblock Arch. Rat. Mech. Anal., \textbf{4} (1960), 273--334.

\bibitem{KR80}
    \newblock E. Kr\"oner,
    \newblock  \emph{Continuum theory of defects},
    \newblock  in ``Physiques des d\'efauts" (ed. R.~Balian),
                Les Houches session XXXV, Course 3, (1980), 219--315.

\bibitem{KR90} 
     \newblock E. Kr\"oner,
     \newblock \emph{The differential geometry of elementary point and line defects in Bravais crystals},
     \newblock  Int. J. Theor. Phys., \textbf{29}, 11, (1990), 1219--1237.

\bibitem{KR92} 
     \newblock E. Kr\"oner,
     \newblock \emph{The internal mechanical state of solids with defects},
     \newblock  Int. J. Solids and Structures, \textbf{29}, 14/15, (1992), 1849--1257.

\bibitem{KR95}
     \newblock E. Kr\"oner,
     \newblock \emph{Dislocations in crystals and in continua: a confrontation},
     \newblock  Int. J. Engng Sci., \textbf{33}, 15, (1995), 2127--2135.

\bibitem{KR96}
     \newblock E. Kr\"oner,
     \newblock \emph{Dislocation theory as a physical field theory},
     \newblock Meccanica, \textbf{31}, (1996), 577--587.

\bibitem{KR01} 
     \newblock E. Kr\"oner,
     \newblock \emph{Benefits and shortcomings of the continuous theory of dislocations},
     \newblock  Int. J. Solids Struc., \textbf{38}, 11, (2001), 1115--1134.

\bibitem{LAZMAUG}
     \newblock M. Lazar and G. Maugin 
     \newblock \textit{Dislocations in gradient elasticity revisited},
    \newblock Proc. R. Soc. A \textbf{462}, (2006), 3465--3480.

\bibitem{Maugin2003}
     \newblock G. Maugin  
	\newblock \emph{Geometry and thermomechanics of structural rearrangements: Ekkehart Kr\"oner's legacy}, 
	  \newblock ZAMM, \textbf{83}, 2, (2003), 75--84.

\bibitem{NOLL67}
     \newblock W. Noll,
     \newblock \emph{Materially uniform bodies with inhomogeneities},
     \newblock Arch. Rational Mech. Anal., \textbf{27}, (1967), 1 -- 32.

\bibitem{NYE53} 
     \newblock J. F. Nye,
       \newblock \emph{Some geometrical relations in dislocated crystals},
           \newblock Acta Metall, \textbf{1}, (1953), 153--162.

\bibitem{PAL2008}
     \newblock M. Palombaro and S. M{\"{u}}ller
     \newblock\emph{Existence of minimizers for a polyconvex energy in a crystal with dislocations},
     \newblock Calc. Var., \textbf{31}, 4, (2008), 473--482.

\bibitem{Phil} 
    \newblock J. Philibert,
    \newblock  ``Atom movements, diffusion and mass transport in solids,"
     \newblock Monographies de physique.,
     \newblock Les \'editions de physique, les Ulis, France., 1988.

\bibitem{Schouten}
    \newblock J. A. Schouten,
    \newblock  ``Ricci-Calculus,"
     \newblock 2nd edn.,
     \newblock Springer-Verlag, Berlin, 1954.

\bibitem{VGetal} 
    \newblock N. {Van Goethem}, A. de Potter, N. Van den Bogaert and F. Dupret,
      \newblock \emph{Dynamic prediction of point defects in {C}zochralski silicon growth. {A}n attempt to reconcile experimental defect diffusion 
coefficients with the ${V/G}$ criterion},
     \newblock J.Phys.Chem.Solids, \textbf{69} (2008), 320--324.

\bibitem{VGD2009}
\newblock N. Van Goethem and F. Dupret,
\newblock \textit{A distributional approach to the geometry of $2{D}$ dislocations at the meso-scale.Parts A and Part B}, 
preprints (2009), arXiv reference: 1003.6021.

\bibitem{VG2010}
\newblock N. Van Goethem,
\newblock \textit{Strain incompatibility in single crystals: Kr\"{o}ner's formula revisited}, 
\newblock J. Elast., DOI: 10.1007/s10659-010-9275-4, (2010).

\bibitem{VG2011}
\newblock N. Van Goethem,
\newblock \textit{A multiscale model for dislocation clusters: from mesoscopic elasticity to macroscopic plasticity}, 
(in preparation), http://ptmat.ptmat.fc.ul.pt/preprints.html.

\bibitem{WANG67}
     \newblock C. C. Wang,
     \newblock \emph{On the geometric structure of simple bodies, a mathematical
foundation for the theory of continuous distributions of dislocations},
     \newblock Arch. Rational Mech. Anal., \textbf{27}, (1967), 33--94.


\end{thebibliography}
\end{document}